\documentclass[a4paper,12pt]{article}

\title{Critical points of a perturbed Otha-Kawasaki functional}
\date{Date}
\author{Matteo Rizzi}
\usepackage[T1]{fontenc}
\usepackage[latin1]{inputenc}
\usepackage[english]{babel}
\usepackage{geometry}
\usepackage{indentfirst} 
\usepackage{graphicx}

\usepackage{amsmath,amsfonts,amssymb,amsthm}
\usepackage{latexsym}

\newcommand{\R}{\mathbb{R}}

\newcommand{\diver}{\text{div}}

\newtheorem{theorem}{Theorem}
\newtheorem{proposition}{Proposition}

\newtheorem{remark}{Remark}

\begin{document}

\maketitle

\begin{abstract}
In the paper, we consider a small perturbation of the Otha-Kawasaki functional and we construct at least four critical points close to suitable translations of the Schwarz P surface with fixed volume. 
\end{abstract}

\tableofcontents

\section{Introduction}
\begin{figure}[htbp]
\centering
\includegraphics[width=150mm]{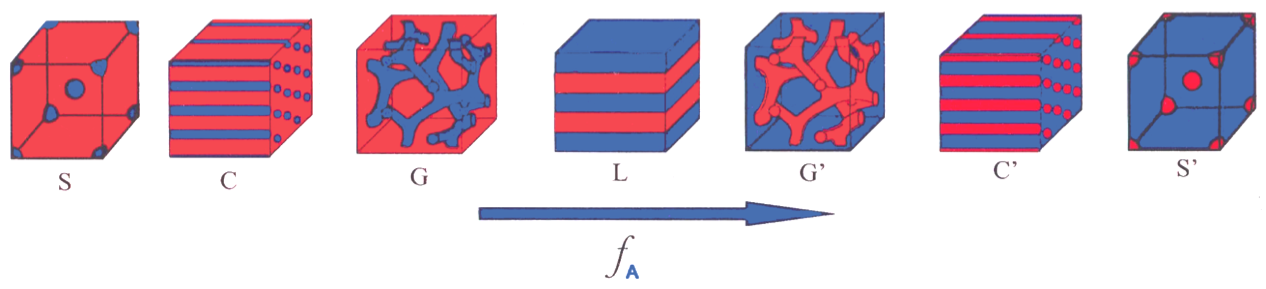}
\caption{The most commonly observed periodic structures are spheres, cylinders, gyroids and lamellae}\label{fig_per}
\end{figure}

A diblock copolymer is a complex molecule where chains of two different kinds of monomers, say A and B, are grafted togheter. Diblock copolymer melts are large collections of diblock copolymers. The experiments show that, above a certain temperature, these melts behave like fluids, that is the monomers are mixed in a disordered way, while below this critical temperature phase separation is observed. Some common periodic structures observed in experiments are spheres, cylinders, gyroids and lamellae (see figure \ref{fig_per}).
These patterns can be found by minimizing some energy. It looks reasonable to describe the phenomenon through an energy given by the sum of the perimeter, that forces the separation surfaces to be minimal, plus some nonlocal term that keeps trace of the long-range interactions between monomers. More explicitly, one can take the functional
\begin{eqnarray}
\mathcal{E}(u):=\frac{1}{2}\int_{\Omega}|\nabla u|dx+\gamma\int_{\Omega}\int_{\Omega}G(x,y)(u(x)-m)(u(y)-m)dxdy
\end{eqnarray}
as an energy. Here $\Omega$ is a bounded domain of $\R^{3}$, that can be seen as the container where the diblock copolimer melt is confined, $u$ is a bounded variation function in $\Omega$ with values in $\{\pm 1\}$ (for instance, we can assume that $u(x)=1$ if there are only monomers of type A at $x$, $u(x)=-1$ if there are only monomers of type B at $x$), $\int_{\Omega}|\nabla u|dx$ is its total variation, or equivalently the perimeter of the set $\{x\in\Omega:u(x)=1\}$, $G$ is the Green's function of $-\Delta$ on $\Omega$, that is the disrtibutional solution to
\begin{eqnarray}\notag
\begin{cases}
-\Delta_{x}G(x,y)=\delta_{y}(x)-\frac{1}{|\Omega|} &\text{in }\Omega\\
\partial_{\nu(x)}G(x,y)=0 &\text{on }\partial\Omega.
\end{cases}
\end{eqnarray}
$G$ turns out to be the sum of the Green's function of $-\Delta$ over $\R^{3}$ and a regular part $R(x,y)$, namely
\begin{eqnarray}\notag
G(x,y)=\frac{c}{|x-y|}+R(x,y),
\end{eqnarray}
(see \cite{RW}). $\gamma\geq 0$ is a parameter depending on the material, that we will assume to be small.\\

This energy appears as the $\Gamma$-limit as $\varepsilon\to 0$ of the approximating functionals
\begin{eqnarray}\notag
\mathcal{E}_{\varepsilon}(u)=\frac{\varepsilon}{2}\int_{\Omega}|\nabla u|^{2}dx+\frac{1}{\varepsilon}\int_{\Omega}\frac{(1-u^{2})^{2}}{4}dx\\\notag
+\frac{16\gamma}{3}\int_{\Omega}\int_{\Omega}G(x,y)(u(x)-m)(u(y)-m)dxdy,
\end{eqnarray}
introduced by Otha and Kawasaki (see \cite{AFM,CM1,CM2,CM3}).\\

In a more geometric way our functional is given by
\begin{eqnarray}
J_{\gamma}(E):=P_{\Omega}(E)+\gamma\int_{\Omega}\int_{\Omega}G(x,y)(u_{E}(x)-m)(u_{E}(y)-m)dxdy
\end{eqnarray}
where
\begin{eqnarray}\notag
E:=\{x\in\Omega:u(x)=1\},
\end{eqnarray}
so that $u_{E}=\chi_{E}-\chi_{\Omega\backslash E}$. The first variation of $J_{\gamma}$ is given by
\begin{eqnarray}
J_{\gamma}^{'}(E)[\varphi]=\int_{\Sigma}(H_{\Sigma}(x)+4\gamma v_{E}(x))\varphi(x)d\sigma(x),
\end{eqnarray}
while its second variation is given by
\begin{eqnarray}
J_{\gamma}^{''}(E)[\varphi]=\int_{\Sigma}L\varphi(x)\varphi(x)d\sigma(x),
\end{eqnarray}
where 
\begin{eqnarray}
L\varphi=-\Delta_{\Sigma}\varphi-|A|^{2}\varphi+8\gamma\int_{\Sigma}G(\cdotp,y)\varphi(y)d\sigma(y)+4\gamma\partial_{\nu}v\varphi.
\end{eqnarray}
Here $\varphi$ is in the space
\begin{eqnarray}
W:=\bigg\{w\in H^{1}(\Sigma):\int_{\Sigma}w(x)\nu_{i}(x)d\sigma(x)=0,&\text{ }1\leq i\leq 3\bigg\},
\label{defW}
\end{eqnarray} 
$\Sigma:=\partial E$ and
\begin{eqnarray}
v_{E}(x):=\int_{T^{3}}G(x,y)(u_{E}(y)-m)dy
\end{eqnarray}
is the unique solution to the problem
\begin{eqnarray}
\begin{cases}
-\Delta v_{E}=u_{E}-m&\text{in $T^{3}$}\\
\int_{T^{3}}v_{E}dx=0.
\end{cases}
\end{eqnarray}
For an explicit computation of the first and the second variation, see for instance \cite{CS}. In the sequel, $\Omega$ will always be the $3$-dimensional torus $T^{3}$, that is the quotient of the cube $[0,1]^{3}$ by the equivalence relation that identifies the opposite faces. It is known that $J_{\gamma}$ is translation invariant, that is $J_{\gamma}(E+\xi)=J_{\gamma}(E)$, for any $\xi\in T^{3}$ (see \cite{AFM},\cite{CS}), thus, once we find a critical point of it, any translation in $T^{3}$ is still critical.\\

There are several results in the literature about critical points of this functional. For instance, an interesting problem is to understand whether all global minimizers are periodic, like the patterns described above (spheres, cylinders, gyroids and lamellae, see Figure \ref{fig_per}). This is known to be true in dimension one (see \cite{M}), but the problem is still open in higher dimension. We refer to \cite{ACO,S} for further results. Some other authors, such as Ren and Wei \cite{RW,RW1,RW2,RW3,RW4}, constructed explicit examples of stable periodic local minimizers, that is with positive second variation. Moreover, Acerbi Fusco and Morini \cite{AFM} showed that any stable critical point is actually a local minimizer with respect to small $L^{1}$ perturbations.

\begin{figure}[htbp]
\centering
\includegraphics[width=45mm]{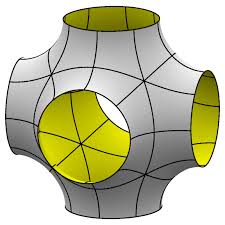}
\caption{Schwarz' P surface}\label{sps}
\end{figure}

Here we add a small linear perturbation that corresponds to an external force $f$ applied to the system, that can be taken to be $C^{0,1}_{loc}(\R^{3})$ and periodic, with triple period $1$. The energy becomes
\begin{eqnarray}
I_{\gamma}(E):=J_{\gamma}(E)+\gamma\int_{\Omega}f(x)u_{E}(x)dx.
\label{defI}
\end{eqnarray}
The additional linear term breakes the translation invariance. We will construct at least four critical points $F_{j}$ of $I_{\gamma}$, $1\leq j\leq 4$, for $\gamma$ small enough, that are close to suitable translations of the Schwarz' P surface $\Sigma$ (see figure \ref{sps}), under the volume constraint
\begin{eqnarray}
\mathcal{L}_{3}(F_{j})=\mathcal{L}_{3}(E),
\end{eqnarray}
where $E$ is the interior of $\Sigma$.
\begin{remark}
The Schwartz P surface can be seen as a periodic surface in $\R^{3}$, with triple period $1$. Moreover, it divides the Torus into two components, an interior and an exterior. In the sequel, $E$ will denote the interior part.
\end{remark}
We will use a technique based on a finite dimensional Lyapunov-Schmidt reduction (see \cite{AM}, Chapter $2.2$), and on the Lusternik-Schnirelman theory (see \cite{NASEP}, Chapter $9$) for the multiplicity.

For $0<\alpha<1$ and for any integer $k\geq 0$, we introduce the H\"{o}lder spaces
\begin{eqnarray}
C^{k,\alpha}_{s}(\overline{\Sigma}):=\{w\in C^{k,\alpha}(\overline{\Sigma}):w(x)=w(T_{j}x),1\leq j\leq 3\},
\end{eqnarray}
where $T_{j}$ are the reflections defined by
\begin{eqnarray}\notag
T_{1}(x_{1},x_{2},x_{3})=(-x_{1},x_{2},x_{3}) & T_{2}(x_{1},x_{2},x_{3})=(x_{1},-x_{2},x_{3}) & T_{3}(x_{1},x_{2},x_{3})=(x_{1},x_{2},-x_{3}). 
\end{eqnarray}
Here it is understood that we have put the origin in the centre of the cube (see Figure \ref{sps}), in such a way that these spaces consist of functions that respect the simmetries of $\Sigma$, that is the simmetries with respect to the coordinate planes $\{x_{j}=0\}$, $1\leq j\leq 3$. We endow these spaces with the norm
\begin{eqnarray}
||w||_{C^{k,\alpha}(\Sigma)}=\sum_{j=0}^{k}||\nabla^{j}w||_{L^{\infty}(\Sigma)}+\sup_{x\neq y}\sup_{|\beta|=k}\frac{|\partial_{\beta}w(x)-\partial_{\beta}w(y)|}{d(x,y)^{\alpha}},
\end{eqnarray}
where $d$ is the geodesics distance on $\Sigma$.
\begin{theorem}
Let $I_{\gamma}$ be defined as in (\ref{defI}) and $\nu(x)$ be the outward-pointing unit normal to the Schwarz P surface $\Sigma$. Then there exists $\gamma_{0}>0$ such that, for any $0<\gamma<\gamma_{0}$, there exist $\xi_{j}\in T^{3}$, $1\leq j\leq 4$, and $w_{\gamma,j}\in C^{2,\alpha}_{s}(\overline{\Sigma})$, with
\begin{eqnarray}
||w_{\gamma,j}||_{C^{2,\alpha}(\Sigma)}\leq c\gamma,
\end{eqnarray}
such that the sets $F_{j}$ defined as the interior of
\begin{eqnarray}
\Gamma_{j}:=\{x+\xi_{j}+\nu(x)w_{\gamma,j}(x):x\in\Sigma\}
\label{formF}
\end{eqnarray}
are critical points of $I_{\gamma}$ under the volume constraint
\begin{eqnarray}
\mathcal{L}_{3}(F_{j})=\mathcal{L}_{3}(E).
\end{eqnarray}
\label{mainth}
\end{theorem}
\begin{remark}
\textit{(i)} If we take $f\equiv 0$, we find a unique critical point $F$, that is the interior of
\begin{eqnarray}
\Gamma:=\{x+\nu(x)w_{\gamma}(x):x\in\partial E\},
\end{eqnarray}
where $w_{\gamma}$ is a small correction, namely $||w_{\gamma}||_{C^{2,\alpha}(\Sigma)}\leq c\gamma$, found by means of the implicit function Theorem (see Remark \ref{ift}). Then any translation $F+\xi$ is still a critical point of $J_{\gamma}$. A similar result was proved by Cristoferi (see \cite{C}, Theorem $4.18$), who constructed a critical point of $J_{\gamma}$ close to any smooth periodic strictly stable constant mean curvature surface.\\

\textit{(ii)} We stated the theorem in the case of $I_{\gamma}$ for simplicity. The same proof should yield existence and multiplicity results also for regular nonlinear perturbations and different coefficients in the nonlocal and forcing terms.
\label{rem_tr_inv}
\end{remark}
A similar result was obtained by Bonacini and Cristoferi \cite{BC}, who studied a nonlocal version of the isoperimetric problem, that is  they considered a small nonlocal perturbation of the perimeter and showed that the unique minimizers $F$ under the volume constraint $\mathcal{L}_{N}(F)=m$ are the balls, provided $m$ is small enough. The critical points we construct here are not necessarily stable, since we apply the Lusternik-Schnirelmann theory (see \cite{NASEP}, chapter $9$).\\   

A crucial tool in the proof is nondegeneracy up to translations of the Jacobi operator of the Schwarz P surface. In \cite{R}, Ross showed that the Schwarz P surface is a critical point of the area and it is \textit{volume preserving stable}, that is it the second variation of the area is non-negative on any normal variation with zero average. More precisely, setting $I_{0}:=P_{\Omega}$, we have
\begin{eqnarray}
I_{0}^{''}(E)(\varphi,\varphi)=\int_{\Sigma}|\nabla_{\Sigma}\varphi|^{2}-|A|^{2}\varphi^{2}d\sigma\geq 0
\end{eqnarray}
for any $\varphi\in H^{1}(\Sigma)$ satisfying
\begin{eqnarray}
\int_{\Sigma}\varphi d\sigma=0,
\end{eqnarray} 
(see Theorem $1$ of \cite{R}). Let $\nu(x)$ denote the exterior unit normal to $\Sigma$ at $x$. Since $I_{0}$ is translation invariant, then $\nu_{i}(x):=(\nu(x),e_{i})$ are Jacobi fields of $\Sigma$, that is they satisfy 
\begin{eqnarray}
-\Delta_{\Sigma}\nu_{i}-|A|^{2}\nu_{i}=0 &\text{in }\Sigma,
\end{eqnarray}
(see \cite{AFM},\cite{CS}). Moreover, Grosse-Brauckmann and Wohlgemuth showed in (\cite{GBW}) that $\Sigma$ is nondegenerate up to translations, that is there are no other nontrivial Jabobi fields. In other words 
\begin{eqnarray}
Ker(I_{0}^{''}(E))=span\{\nu_{i}\}_{1\leq i\leq 3}.
\label{ker_Jacobi}
\end{eqnarray}
\begin{remark}
Let us observe that the $\nu_{i}$'s are linearly independent. In fact, if not, there would exist a constant vector $b=(b_{1},b_{2},b_{3})\neq 0$ such that $0=(b,\nu(x))$ for any $x\in\Sigma$, but this contradicts the geometry of $\Sigma$.
\label{nu_i_li}
\end{remark}
We note that the $\nu_{i}$'s have zero average, since
\begin{eqnarray}
\int_{\Sigma}\nu_{i}(x)d\sigma(x)=\int_{T^{3}}\diver e_{i}=0.
\end{eqnarray}
In addition, we decompose $H^{1}(\Sigma)$ into the orthogonal sum 
\begin{eqnarray}
H^{1}(\Sigma)=span\{\nu_{i}\}_{1\leq i\leq 3}+W,
\end{eqnarray}
(see (\ref{defW}) for the definition of $W$), and we define
\begin{eqnarray}
W^{0}:=\bigg\{w\in W:\int_{\Sigma}w(x)d\sigma(x)=0\bigg\}.
\end{eqnarray}
The above discussion can be rephrased by saying that
\begin{eqnarray}
\int_{\Sigma}|\nabla_{\Sigma}w|^{2}-|A|^{2}w^{2}d\sigma\geq c||w||_{H^{1}(\Sigma)}^{2} &\text{for any }w\in W^{0}. 
\label{Jacobi_pos}
\end{eqnarray}

\textbf{Aknowledgments} The author is supported by the PRIN project \textit{Variational and
perturbative aspects of nonlinear differential problems}. The author is also particularly
grateful to F. Mahmoudi for his precious collaboration.

\section{The proof of Theorem \ref{mainth}: Lyapunov-Schmidt reduction}
We need to find at least four sets $F$ of the form (\ref{formF}) and a Lagrange multiplier $\lambda\in\R$ such that
\begin{eqnarray}
H_{\partial F}(y)+4\gamma v_{F}(y)+\gamma f(y)=\lambda &\text{ }\forall y\in\partial F,
\label{maineq}
\end{eqnarray}
or equivalently
\begin{eqnarray}
I^{'}_{\gamma}(F)=\lambda.
\end{eqnarray}
Exploiting the variational nature of the problem and the fact that $H_{\Sigma}=0$, equation (\ref{maineq}) is equivalent to
\begin{eqnarray}
\lambda=4\gamma v_{E}(x)+Lw(x)+Q(w)(x)+\gamma f(y), &\forall x\in\Sigma,
\label{maineq_taylor}
\end{eqnarray}
where $y$ is seen as a function of $x$ depending on the parameter $\xi$, namely $y=x+\xi+w(x)\nu(x)$, and
\begin{eqnarray}
Q(w):=J^{'}_{\gamma}(F)-J^{'}_{\gamma}(E)-J^{''}_{\gamma}(E)w.
\end{eqnarray}
Writing
\begin{eqnarray}
Lw=-\Delta_{\Sigma}w-|A|^{2}w+\gamma\tilde{L}w,
\end{eqnarray}
where
\begin{eqnarray}
\tilde{L}w=8\int_{\Sigma}G(\cdotp,\zeta)w(\zeta)d\sigma(\zeta)+4\partial_{\nu}v_{E}w,
\end{eqnarray}
we can see that (\ref{maineq_taylor}) is equivalent to
\begin{eqnarray}
-\Delta_{\Sigma}w-|A|^{2}w=\lambda+\mathcal{F}(\gamma,\xi,w),
\label{maineq_w}
\end{eqnarray}
where the nonlinear functional $\mathcal{F}$ is given by
\begin{eqnarray}
\mathcal{F}(\gamma,\xi,w)(x)=-4\gamma v_{E}(x)-\gamma\tilde{L}w(x)-Q(w)(x)-\gamma f(y), &\forall x\in\Sigma.
\end{eqnarray}
The unknowns are the function $w$, $\xi\in T^{3}$ and $\lambda\in\R$.

\subsection{The volume constraint}
Now we will consider the relation between the volume of $F$ and $w$. In order to do so, we point out that there exists a global parametrization 
\begin{eqnarray}
\phi:Y\to\Sigma,
\label{def_par}
\end{eqnarray}
defined on an open set $Y\in\R^{2}$ (see \cite{GK}, section $3$), that induces a change of coordinates on a neighbourhood of $\Sigma$ given by
\begin{eqnarray}
X(\text{y}_{1},\text{y}_{2},z):=\phi(\text{y}_{1},\text{y}_{2})+z\nu(\text{y}_{1},\text{y}_{2}),
\end{eqnarray}
where, with an abuse of notation, $\nu(\text{y}_{1},\text{y}_{2})$ is the outward-pointing unit normal to $\Sigma$ at $\phi(\text{y}_{1},\text{y}_{2})$. The volume of $F$ is given by
\begin{eqnarray}\notag
\mathcal{L}_{3}(F)=\mathcal{L}_{3}(E)+\int_{Y}d\text{y}\int_{0}^{w(\text{y})}\det JX(\text{y},z)dz,
\end{eqnarray}
where $JX$ is the Jacobian of $X$. We expand
\begin{eqnarray}\notag
\det JX(\text{y},z)=\det JX(\text{y},0)+zA(\text{y})+z^{2}B(\text{y}),
\end{eqnarray}
thus we get
\begin{eqnarray}\notag
\mathcal{L}_{3}(F)=\mathcal{L}_{3}(E)+\int_{Y}d\text{y}\int_{0}^{w(\text{y})}\bigg(\det JX(\text{y},0)+zA(\text{y})+z^{2}B(\text{y})\bigg)dz\\\notag
=\mathcal{L}_{3}(E)+\int_{Y}\det JX(\text{y},0)w(\text{y})d\text{y}+\int_{Y}\bigg(\frac{1}{2}w^{2}(\text{y})A(\text{y})+\frac{1}{3}w^{3}(\text{y})B(\text{y})\bigg)d\text{y}.
\end{eqnarray}
Since $\det JX(\text{y},0)=(\nu(\text{y}),\partial_{\text{y}_{1}}\phi\times\partial_{\text{y}_{2}}\phi)\neq 0$ for any y$\in Y$, 
\begin{eqnarray}
\mathcal{L}_{3}(F)=\mathcal{L}_{3}(E)+\int_{\Sigma}w(x)d\sigma(x)+\int_{\Sigma}\tilde{Q}(x,w(x))d\sigma(x),
\end{eqnarray}
where
\begin{eqnarray}
\tilde{Q}(x,w)=\frac{1}{\det JX(x)}\bigg(\frac{1}{2}w^{2}(x)A(x)+\frac{1}{3}w^{3}(x)B(x)\bigg).
\end{eqnarray}
Therefore the volume constraint is equivalent to an equation of the form
\begin{eqnarray}
\int_{\Sigma}w(x)dx=-\int_{\Sigma}\tilde{Q}(x,w(x))d\sigma(x).
\label{constraint}
\end{eqnarray}

\subsection{The auxiliary equation}
The aim is to solve (\ref{maineq_w}) under the volume constraint (\ref{constraint}). However, since, by (\ref{ker_Jacobi}) and (\ref{Jacobi_pos}), the Jacobi operator $-\Delta _{\Sigma}-|A|^{2}$ is non degenerate up to translations, we can actually solve the system
\begin{eqnarray}
\begin{cases}
-\Delta_{\Sigma}w-|A|^{2}w=\lambda+P\mathcal{F}(\gamma,\xi,w) &\text{in }\Sigma\\
\partial_{n} w=0 &\text{on }\partial\Sigma,\\
\int_{\Sigma}wd\sigma=-\int_{\Sigma}\tilde{Q}(x,w(x))d\sigma(x),
\end{cases}
\label{sist_aux}
\end{eqnarray}
where $P:L^{2}(\Sigma)\to\tilde{W}$ is the projection onto the space
\begin{eqnarray}
\tilde{W}:=\bigg\{\varphi\in L^{2}(\Sigma):\int_{\Sigma}\varphi(x)\nu_{i}(x)d\sigma(x)=0,&1\leq i\leq 3\bigg\},
\end{eqnarray}
$\partial_{n} w:=(\nabla_{\Sigma}w,n)$ and $n$ is the outward pointing unit normal to $\partial\Sigma$ in $\Sigma$. This will be done by a fixed point argument in the following Proposition, proved in section $3$. 
\begin{proposition}
For any $\xi\in T^{3}$ and for any $\gamma$ sufficiently small, there exists a unique solution $(w_{\gamma,\xi},\lambda_{\gamma,\xi})\in C^{2,\alpha}_{s}(\overline{\Sigma})\times\R$ to problem (\ref{sist_aux}) satisfying 
\begin{eqnarray}
||w_{\gamma,\xi}||_{C^{2,\alpha}(\Sigma)}+|\lambda_{\gamma,\xi}|\leq C\gamma,\\
\int_{\Sigma}w(x)\nu_{i}(x)d\sigma(x)=0,&\text{ }1\leq i\leq 3,
\end{eqnarray}
for some constant $C>0$. Moreover, the solution is Lipschitz continuous with respect to the parameter $\xi$, that is
\begin{eqnarray}
||w_{\gamma,\xi_{1}}-w_{\gamma,\xi_{2}}||_{C^{2,\alpha}(\Sigma)}+|\lambda_{\gamma,\xi_{1}}-\lambda_{\gamma,\xi_{1}}|\leq C\gamma|\xi_{1}-\xi_{2}|,& \forall\xi_{1},\xi_{2}\in T^{3}.
\label{lip_xi}
\end{eqnarray}
\label{prop_aux}
\end{proposition}
\begin{remark}
If we take $f\equiv 0$, in order to get the right correction $w$, we just solve (\ref{sist_aux}) for $\xi=0$, due to the translation invariance of $J_{\gamma}$ (see Remark \ref{rem_tr_inv}). We do not need the Lyapunov-Schmidt reduction.
\label{ift}
\end{remark}

\subsection{The bifurcation equation}
In order to conclude the proof of Theorem \ref{mainth}, we have to find at least four points $\xi\in T^{3}$ such that 
$(Id-P)\mathcal{F}(\gamma,\xi,w_{\gamma,\xi})(x)=0$, or equivalently
\begin{eqnarray}
\int_{\Sigma}\mathcal{F}(\gamma,\xi,w_{\gamma,\xi})(x)\nu_{i}(x)d\sigma(x)=0,
\label{bifo_eq}
\end{eqnarray}
for $i=1,2,3$.\\ 

Since $\partial_{n}w_{\gamma,\xi}=0$ on $\partial\Sigma$ and the same is true for the $\nu_{i}$'s, an integration by parts yields
\begin{eqnarray}\notag
\int_{\Sigma}(-\Delta_{\Sigma}w_{\gamma,\xi}-|A|^{2}w_{\gamma,\xi})(x)\nu_{i}(x)d\sigma(x)=0,
\end{eqnarray}
for $i=1,2,3$, thus by (\ref{sist_aux}) we can see that $w$ solves 
\begin{eqnarray}
P(I_{\gamma}^{'}(F)-\lambda)=0,
\label{aux_eq}
\end{eqnarray}
or equivalently
\begin{eqnarray}
I_{\gamma}^{'}(F)-\lambda=\sum_{i=1}^{3}A_{i,\gamma,\xi}\nu_{i}.
\label{aux_eq_Ai}
\end{eqnarray}
Since, by construction,
\begin{eqnarray}\notag
\mathcal{F}(\gamma,\xi,w_{\gamma,\xi})=-I_{\gamma}^{'}(F)-\Delta_{\Sigma}w_{\gamma,\xi}-|A|^{2}w_{\gamma,\xi}=\\\notag
-I_{\gamma}^{'}(F)+\lambda+P\mathcal{F}(\gamma,\xi,w)=-\sum_{i=1}^{3}A_{i,\gamma,\xi}\nu_{i}+P\mathcal{F}(\gamma,\xi,w),
\end{eqnarray}
and (\ref{aux_eq_Ai}) holds, we can see that (\ref{bifo_eq}) is equivalent to
\begin{eqnarray}
A_{i,\gamma,\xi}=0 &\text{for }i=1,2,3.
\label{Ai=0}
\end{eqnarray}
Equation (\ref{Ai=0}) is solvable thanks to the Lusternik-Schnirelmann theory and the compactness of the Torus. We recall that the Torus $T^{3}$ has category $4$ (see \cite{NASEP}, example $9.4$, \textit{(iii)}).
\begin{proposition}
Equation (\ref{Ai=0}) is satisfied if $\xi$ is a critical point of the function $\Phi_{\gamma}:T^{3}\to\R$ defined by
\begin{eqnarray}
\Phi_{\gamma}(\xi):=I_{\gamma}(F),
\label{def_phi}
\end{eqnarray}
where $F$ is the interior of
\begin{eqnarray}\notag
\Gamma:=\{x+\xi+w_{\gamma,\xi}(x)\nu(x):x\in\Sigma\}.
\end{eqnarray}
\label{prop_bifo}
\end{proposition}
The proof of Proposition \ref{prop_bifo} will be carried out in Section $4$. It is possible to see that $\Phi_{\gamma}$ actually admits at least $4$ critical points, due to Theorem $9.10$ of \cite{NASEP} applied to $I_{\gamma}$, with $M=T^{3}$. The compactness of the torus $T^{3}$ is crucial, since it guarantees that $I_{\gamma}$ is bounded from below on $M$ and the Palais-Smale condition is satisfied.

\section{Solving the auxiliary equation}
The aim of this section is to prove Proposition \ref{prop_aux}. First, in Section $3,1$, we will treat the corresponding linear problem, then, in Section $3,2$, we will solve problem (\ref{sist_aux}) by a fixed point argument.
\subsection{The linear problem}
\begin{proposition}
Let $a\in\R$ and $\varphi\in C^{0,\alpha}_{s}(\overline{\Sigma})$ be such that
\begin{eqnarray}
\int_{\Sigma}\varphi\nu_{i}d\sigma=0 &\text{for }i=1,2,3.
\label{ort_varphi}
\end{eqnarray}
Then there exists a unique solution $(w,\lambda)=\Psi(\varphi,a)\in C^{2,\alpha}_{s}(\overline{\Sigma})\times\R$ to the problem
\begin{eqnarray}
\begin{cases}
-\Delta_{\Sigma}w-|A|^{2}w=\lambda+\varphi&\text{in }\Sigma\\
\partial_{n}w=0 &\text{on }\partial\Sigma\\
\int_{\Sigma}w\nu_{i}d\sigma=0 &\text{for }1\leq i\leq 3,\\
\int_{\Sigma}wd\sigma=a.
\end{cases}
\label{lin_pb}
\end{eqnarray}
Moreover, we have the stimate
\begin{eqnarray}
||w||_{C^{2,\alpha}(\Sigma)}+|\lambda|\leq c(||\varphi||_{C^{0,\alpha}(\Sigma)}+|a|).
\end{eqnarray}
\label{prop_lin}
\end{proposition}
\begin{remark}
Since the $\nu_{i}$'s are linearly independent (see Remark \ref{nu_i_li}), then the matrix
\begin{eqnarray}
L_{ki}:=\int_{\Sigma}\nu_{k}\nu_{i}d\sigma
\label{defL}
\end{eqnarray}
is invertible (for a detailed proof, see the appendix).
\label{Lki_invt}
\end{remark}
\begin{proof}
\textit{Step (i): existence and uniqueness.}\\ 
First we look for a weak solution $w\in W$. We write any $w\in W$ as
\begin{eqnarray}\notag
w=w_{0}+\frac{1}{|\Sigma|}\int_{\Sigma}wd\sigma,
\end{eqnarray}
with $w_{0}\in\ W^{0}$. The linear problem can be rephrased as follows
\begin{eqnarray}
\begin{cases}
-\Delta_{\Sigma}w_{0}-|A|^{2}w_{0}=\lambda+\varphi+|A|^{2}\frac{a}{|\Sigma|} &\text{in }\Sigma\\
\int_{\Sigma}w_{0}=0.
\end{cases}
\label{lin_pb_0}
\end{eqnarray}
We note that the right-hand side of (\ref{lin_pb_0}) is orthogonal to $\nu_{i}$, for $i=1,2,3$, due to the fact that
\begin{eqnarray}
\int_{\Sigma}|A|^{2}\frac{a}{|\Sigma|}\nu_{i}(x)d\sigma=
\int_{\Sigma}\bigg(\Delta_{\Sigma}\nu_{i}+|A|^{2}\nu_{i}\bigg)(x)\frac{a}{|\Sigma|}d\sigma=0,
\label{Delta_nu_i_0}
\end{eqnarray}
since $\partial_{n}\nu_{i}=0$ on $\partial\Sigma$, and
\begin{eqnarray}
\int_{\Sigma}\nu_{i}(x)d\sigma(x)=\int_{E}\diver(e_{1})dx=0.
\label{av_nui}
\end{eqnarray}
In addition, the norm defined by
\begin{eqnarray}
||w||=\int_{\Sigma}|\nabla_{\Sigma}w|^{2}-|A|^{2}w^{2}
\end{eqnarray}
is equivalent to the $H^{1}(\Sigma)$-norm on $W^{0}$, thus the functional
\begin{eqnarray}\notag
I(w)=\int_{\Sigma}|\nabla_{\Sigma}w|^{2}-|A|^{2}w^{2}d\sigma-\int_{\Sigma}\bigg(\varphi+|A|^{2}\frac{a}{|\Sigma|}\bigg)wd\sigma
\end{eqnarray}
is bounded from below by
\begin{eqnarray}
I(w)\geq c||w||_{H^{1}(\Sigma)}^{2}-||\varphi||_{L^{2}(\Sigma)}||w||_{H^{1}(\Sigma)},
\end{eqnarray}
on $W^{0}$, hence it is coercive on it. Moreover, this functional is also w.l.s.c. and strictly convex on $W^{0}$, therefore any minimizing sequence $w_{k}\in W^{0}$ weakly converges, up to subsequence, to the unique minimizer $w_{0}\in W^{0}$, which satisfies the Euler-Lagrange equation   
\begin{eqnarray}\notag
\int_{\partial\Sigma}\partial_{n}wvd\sigma_{\partial\Sigma}+\int_{\Sigma}(-\Delta_{\Sigma}w_{0}-|A|^{2}w_{0})vd\sigma=\\
\lambda\int_{\Sigma}vd\sigma+\sum_{i=1}^{3}\beta_{i}\int_{\Sigma}\nu_{i}vd\sigma+\int_{\Sigma}\varphi vd\sigma+\int_{\Sigma}|A|^{2}\frac{a}{|\Sigma|}vd\sigma,
\label{weak_sol}
\end{eqnarray}
for any $v\in H^{1}(\Sigma)$, for some Lagrange multipliers $\lambda,\beta_{i}\in\R$. Since $\varphi\in C^{0,\alpha}(\overline{\Sigma})$, then $w\in C^{2,\alpha}(\overline{\Sigma})$ (see for instance \cite{N}). Taking the test functions $v\in C^{1}_{c}(\Sigma)$, we can see that $w$ satsfies 
\begin{eqnarray}\notag
-\Delta_{\Sigma}w_{0}-|A|^{2}w_{0}=\lambda+\sum_{i=1}^{3}\beta_{i}\nu_{i}+\varphi+|A|^{2}\frac{a}{|\Sigma|}&\text{in }\Sigma,
\end{eqnarray}
in the classical sense. Taking now $v\in C^{1}(\overline{\Sigma})$, we can see that the Neumann boundary condition is satisfied in the classical sense too. Moreover, $w$ respects the required simmetries because of the symmetries of the laplacian and uniqueness. Taking $\nu_{j}$ as a test function in (\ref{weak_sol}), using (\ref{av_nui}), (\ref{ort_varphi}), (\ref{Delta_nu_i_0}) the Neumann boundary condition and the fact that $\partial_{n}\nu_{i}=0$ on $\partial\Sigma$, we get 
\begin{eqnarray}\notag
\sum_{i=1}^{3}\beta_{i}\int_{\Sigma}\nu_{i}\nu_{j}d\sigma=0,
\end{eqnarray}
therefore by Remark \ref{Lki_invt}, $\beta_{i}=0$.\\

\textit{Step (ii): Regularity estimates.}\\
Multiplying (\ref{lin_pb_0}) by $w_{0}$, integrating by parts and using (\ref{Jacobi_pos}), the Neumann boundary conditions and H\"{o}lder's inequality, we can see that
\begin{eqnarray}\notag
c||w_{0}||_{H^{1}(\Sigma)}^{2}\leq\int_{\Sigma}|\nabla_{\Sigma}w_{0}|^{2}-|A|^{2}w_{0}^{2}d\sigma=\int_{\Sigma}\varphi w_{0}d\sigma+\frac{a}{|\Sigma|}\int_{\Sigma}|A|^{2}w_{0}d\sigma\leq\\\notag
||w_{0}||_{L^{2}(\Sigma)}(||\varphi||_{L^{2}(\Sigma)}+\tilde{c}|a|)\leq||w_{0}||_{H^{1}(\Sigma)}(||\varphi||_{L^{2}(\Sigma)}+\tilde{c}|a|).
\end{eqnarray}
Since $||w||_{H^{1}(\Sigma)}^{2}=||w_{0}||_{H^{1}(\Sigma)}^{2}+a^{2}$, then
\begin{eqnarray}\notag
||w||_{H^{1}(\Sigma)}\leq c(||\varphi||_{L^{2}(\Sigma)}+|a|).
\end{eqnarray}
In order to estimate $\lambda$, we integrate (\ref{lin_pb}) and we get
\begin{eqnarray}\notag
\lambda|\Sigma|+\int_{\Sigma}\varphi d\sigma=-\int_{\Sigma}|A|^{2}wd\sigma,
\end{eqnarray}
since, by the Neumann boundary conditions,
\begin{eqnarray}
\int_{\Sigma}\Delta_{\Sigma}wd\sigma=\int_{\partial\Sigma}\partial_{n}wd\sigma_{\partial\Sigma}=0,
\end{eqnarray}
thus
\begin{eqnarray}\notag
|\lambda|\leq c(||\varphi||_{L^{2}(\Sigma)}+||w||_{L^{2}(\Sigma)}).
\end{eqnarray}
To sum up, we have the estimate
\begin{eqnarray}
|\lambda|+||w||_{H^{1}(\Sigma)}\leq c(||\varphi||_{L^{2}(\Sigma)}+|a|),
\label{contH1}
\end{eqnarray}
In order to get the estimate with respect to the norms we are interested in, we point out that, by the Sobolev embeddings
\begin{eqnarray}\notag
||w||_{L^{\infty}(B_{\delta}(x))}\leq c||w||_{W^{2,2}(B_{\delta}(x))}\leq c(||w||_{L^{2}(B_{2\delta}(x))}+||\varphi+\lambda||_{L^{2}(B_{2\delta}(x))}+|a|)\leq\\\notag
c(||\varphi||_{L^{2}(\Sigma)}+|a|),
\end{eqnarray}
for any $\delta>0$ small but fixed and $x\in\Sigma$ such that $d(x,\partial\Sigma)>\delta$ (here, $B_{\delta}(x)$ is the geodesic ball of radius $\delta$ centered at $x$ in $\Sigma$). In particular, 
\begin{eqnarray}\notag
||w||_{L^{\infty}(\Sigma)}\leq c(||\varphi||_{L^{\infty}(\Sigma)}+|a|).
\end{eqnarray}
By the H\"{o}lder's regularity estimates, we conclude that,
\begin{eqnarray}\notag
||w||_{C^{2,\alpha}(\Sigma)}\leq c(||w||_{L^{\infty}(\Sigma)}+||\varphi+\lambda||_{C^{0,\alpha}(\Sigma)})\leq c(||\varphi||_{C^{0,\alpha}(\Sigma)}+|a|),
\end{eqnarray}
(see \cite{GT}, Chapter $6$, Theorem $6.30$). Since the same is true for $|\lambda|$, the proof is over.
\end{proof}

\subsection{The proof of Proposition \ref{prop_aux}: a fixed point argument}
Now we are ready to show existence, uniqueness and Lipschitz continuity with respect to $\xi$ of the solution $(w,\lambda)$ to (\ref{sist_aux}).\\

\textit{Step (i): Existence and uniqueness.}\\

We solve our problem by a fixed point argument. In fact the map
\begin{eqnarray}\notag
T(w,\lambda)=\Psi\bigg(P\mathcal{F}(\gamma,\xi,w),-\int_{\Sigma}\tilde{Q}(w)\bigg)
\end{eqnarray}
is a contraction on the product $B\times\Lambda$, where $\Lambda=(-C\gamma,C\gamma)$ and
\begin{eqnarray}
B:=\{w\in W\cap C^{2,\alpha}_{s}(\overline{\Sigma}):\partial_{n}w=0\text{ on $\partial\Sigma$, }||w||_{C^{2,\alpha}(\Sigma)}<C\gamma\},
\end{eqnarray}
provided $C$ is large enough. In fact
\begin{eqnarray}\notag
||\mathcal{F}(\gamma,\xi,w)||_{C^{0,\alpha}(\Sigma)}\leq \gamma(4||v_{E}||_{C^{2,\alpha}(\Sigma)}+||f||_{C^{0,\alpha}(\Sigma)})+c\gamma||w||_{C^{2,\alpha}(\Sigma)}\leq\\\notag \gamma(4||v_{E}||_{C^{2,\alpha}(\Sigma)}+||f||_{C^{0,\alpha}(\Sigma)})+cC\gamma^{2}<C\gamma
\end{eqnarray}
provided $C>2(4||v_{E}||_{C^{2,\alpha}(\Sigma)}+||f||_{C^{0,\alpha}(\Sigma)})$ and $\gamma$ is small enough. Similarly, we can see that $\mathcal{F}(\gamma,\xi,w)$ is Lipschitz continuous in $w$ with Lipschitz constant of order $\gamma$.

In addition, the second component fulfills
\begin{eqnarray}\notag
\bigg|\int_{\Sigma}\tilde{Q}(w)\bigg|\leq c||w||_{C^{2,\alpha}(\Sigma)}^{2}\leq cC^{2}\gamma^{2}<C\gamma
\end{eqnarray}
if $\gamma$ is small enough, and the same is true for the Lipschitz constant.\\

\textit{Lipschitz continuity with respect to $\xi$.}\\

In order to prove (\ref{lip_xi}), we point out that, if we set $w_{i}:=w_{\xi_{i}}$ and $y_{i}:=x+\xi_{i}+w_{i}(x)\nu(x)$, for $i=1,2$,
\begin{eqnarray}\notag
||f(y_{1})-f(y_{2})||_{C^{0,\alpha}(\Sigma)}\leq c(|\xi_{1}-\xi_{2}|+||w_{1}-w_{2}||_{C^{2,\alpha}(\Sigma)})\\\notag
||\tilde{L}w_{1}-\tilde{L}w_{2}||_{C^{0,\alpha}(\Sigma)}\leq c||w_{1}-w_{2}||_{C^{2,\alpha}(\Sigma)}
\end{eqnarray}
and
\begin{eqnarray}\notag
||Q(w_{1})-Q(w_{2})||_{C^{0,\alpha}(\Sigma)}\leq c(||w_{1}||_{C^{2,\alpha}(\Sigma)}+||w_{2}||_{C^{2,\alpha}(\Sigma)})||w_{1}-w_{2}||_{C^{2,\alpha}(\Sigma)}\leq\\\notag cC\gamma||w_{1}-w_{2}||_{C^{2,\alpha}(\Sigma)}.
\end{eqnarray}
Similarly, we can show that
\begin{eqnarray}\notag
\bigg|\int_{\Sigma}(\tilde{Q}(w_{1})-\tilde{Q}(w_{2}))d\sigma\bigg|\leq c\gamma||w_{1}-w_{2}||_{C^{2,\alpha}(\Sigma)},
\end{eqnarray}
thus, applying $\Psi$,
\begin{eqnarray}\notag
|\lambda_{1}-\lambda_{2}|+||w_{1}-w_{2}||_{C^{2,\alpha}(\Sigma)}\leq c\gamma(||w_{1}-w_{2}||_{C^{2,\alpha}(\Sigma)}+|\xi_{1}-\xi_{2}|).
\end{eqnarray}
In conclusion, for $\gamma$ small enough,
\begin{eqnarray}\notag
|\lambda_{1}-\lambda_{2}|+\frac{1}{2}||w_{1}-w_{2}||_{C^{2,\alpha}(\Sigma)}\leq c\gamma|\xi_{1}-\xi_{2}|.
\end{eqnarray}

\section{Solving the bifurcation equation.}
The parametrization $\phi:Y\to\Sigma$ of $\Sigma$ introduced in (\ref{def_par}) induces a parmetrization $\beta:Y\to\Gamma:=\partial F$ given by
\begin{eqnarray}
\beta(\text{y}_{1},\text{y}_{2}):=\phi(\text{y}_{1},\text{y}_{2})+\xi+w_{\gamma,\xi}(\text{y}_{1},\text{y}_{2})\nu(\text{y}_{1},\text{y}_{2}).
\end{eqnarray}
The volume element can be expressed in terms of $\phi$ in this way
\begin{eqnarray}\notag
|\beta_{\text{y}_{1}}\times\beta_{\text{y}_{2}}|=|\phi_{\text{y}_{1}}\times\phi_{\text{y}_{2}}|+L^{1}_{\xi}w_{\xi}+Q^{1}_{\xi}w_{\xi},
\label{exp_dvol}
\end{eqnarray}
where $L^{1}_{\xi}$ depends linearly on $w_{\xi}$ and on its gradient and $Q^{1}_{\xi}$ is quadratic in the same quantites. More precisely, they satisfy the estimates
\begin{eqnarray}
\begin{cases}
|L^{1}_{\xi}w|\leq c||w||_{C^{2,\alpha}(\Sigma)}\\
|Q^{1}_{\xi}(w)|\leq c||w||_{C^{2,\alpha}(\Sigma)}^{2}.
\end{cases}
\label{remainder}
\end{eqnarray}
Using the Taylor expansion of the function $\frac{1}{1+s}$, we can show that the outward-pointing unit normal to $\Gamma$ is  
\begin{eqnarray}
\nu_{\Gamma}=\frac{\beta_{\text{y}_{1}}\times\beta_{\text{y}_{2}}}{|\beta_{\text{y}_{1}}\times\beta_{\text{y}_{2}}|}
=\frac{\phi_{\text{y}_{1}}\times\phi_{\text{y}_{2}}}{|\phi_{\text{y}_{1}}\times\phi_{\text{y}_{2}}|}
+\tilde{L}^{1}_{\xi}w_{\xi}+\tilde{Q}^{1}_{\xi}w_{\xi}=\\\notag
\nu+\tilde{L}^{1}_{\xi}w_{\xi}+\tilde{Q}^{1}_{\xi}w_{\xi},
\label{exp_nu}
\end{eqnarray}
with $\tilde{L}^{1}_{\xi}$ and $\tilde{Q}^{1}_{\xi}$ satisfying (\ref{remainder}).

Now we point out that, if $\xi$ is a critical point of $\Phi_{\gamma}$, then 
\begin{eqnarray}
\partial_{\xi_{i}}\Phi_{\gamma}(\xi)=0.
\label{crit_Phi}
\end{eqnarray}
We will rephrase this fact in a more convenient way, that will be more suitable for the forthcoming computations. We define the one-parameter family of diffeomorphisms
\begin{eqnarray}\notag
y_{t}:Y\to\R^{3}
\end{eqnarray}
by 
\begin{eqnarray}
y_{t}(\text{y}_{1},\text{y}_{2}):=\phi(\text{y}_{1},\text{y}_{2})+\xi+te_{i}+w_{\gamma,\xi+te_{i}}(\text{y}_{1},\text{y}_{2})
\nu(\text{y}_{1},\text{y}_{2}),
\end{eqnarray}
for $i=1,2,3$; $\Gamma_{t}:=y_{t}(Y)$ is the image of $y_{t}$. By construction, $\Gamma_{t}$ is actually a submanifold of $T^{3}$ and $\Gamma_{0}=\Gamma$. In terms of $\Gamma_{t}$, condition (\ref{crit_Phi}) is equivalent to
\begin{eqnarray}
\frac{d}{dt}I_{\gamma}(\Gamma_{t})|_{t=0}=0.
\end{eqnarray}
By a result of Fall and Mahmoudi (see \cite{FM}), 
\begin{eqnarray}
0=\frac{d}{dt}I_{\gamma}(\Gamma_{t})|_{t=0}=\int_{\Gamma}(H_{\Gamma}+4\gamma v_{F}+f)(\zeta,\nu_{\Gamma})d\sigma_{\Gamma}+\frac{1}{|\partial\Gamma|}\int_{\partial\Gamma}(\zeta,\nu^{\Gamma}_{\partial\Gamma})ds,
\end{eqnarray}
where 
\begin{eqnarray}
\zeta=\frac{d}{dt}y_{t}(x)|_{t=0}=e_{i}+\partial_{\xi_{i}}w_{\xi}\nu.
\end{eqnarray}
and $\nu^{\Gamma}_{\partial\Gamma}$ is the unit normal to $\partial\Gamma$ in $\Gamma$. The boundary term vanishes by periodicity and by the symmetries of the problem. Using the parametrization $\beta$ of $\Gamma$ and expansions (\ref{exp_nu}) and (\ref{exp_dvol}), the latter relation becomes
\begin{eqnarray}\notag
\int_{Y}\bigg\{(H_{\Gamma}+4\gamma v_{F}+f)(\beta(\text{y}_{1},\text{y}_{2}))\\\notag
(e_{i}+\partial_{\xi_{i}}w_{\xi}\nu,\nu+\tilde{L}^{1}_{\xi}w_{\xi}+\tilde{Q}^{1}_{\xi}w_{\xi})\\\notag
(|\phi_{x}\times\phi_{y}|+L^{1}_{\xi}w_{\xi}+Q^{1}_{\xi}w_{\xi})\bigg\}d\text{y}_{1}d\text{y}_{2}=0.
\end{eqnarray}
By the auxiliary equation, we know that
\begin{eqnarray}
(H_{\Gamma}+4\gamma v_{F}+f)(\beta(\text{y}_{1},\text{y}_{2}))=\sum_{k=1}^{3}A_{k,\gamma,\xi}\nu_{k}(\text{y}_{1},\text{y}_{2})+\lambda,
\end{eqnarray}
thus
\begin{eqnarray}
\sum_{k=1}^{3}A_{k,\gamma,\xi}\bigg(\int_{\Sigma}\nu_{k}\nu_{i}d\sigma+b_{ki}\bigg)+\lambda\int_{\Gamma}(\zeta,\nu_{\Gamma})d\sigma_{\Gamma}=0, &\text{for }i=1,2,3,
\end{eqnarray}
with $b_{ki}=O(\gamma)$. Moreover, once again by \cite{FM}, we know that
\begin{eqnarray}\notag
\frac{d}{dt}\mathcal{L}_{3}(F_{t})=\int_{\Gamma}(\zeta,\nu_{F})d\sigma_{\Gamma},
\end{eqnarray}
hence, by the volume constraint, 
\begin{eqnarray}\notag
\int_{\Gamma}(\zeta,\nu_{F})d\sigma_{\Gamma}=0,
\end{eqnarray}
thus we get
\begin{eqnarray}
\sum_{k=1}^{3}A_{k,\gamma,\xi}\bigg(\int_{\Sigma}\nu_{k}\nu_{i}d\sigma+b_{ki}\bigg)=0, &\text{for }i=1,2,3.
\end{eqnarray}
Since the matrix $L_{ki}$ is invertible (see Remark \ref{Lki_invt}) and the coefficients $b_{ki}$ are small, the matrix $L_{ki}+b_{ki}$ is invertible too, therefore $A_{k,\gamma,\xi}=0$ for $k=1,2,3$.

\section{Appendix}
\textit{Proof of Remark \ref{Lki_invt}}\\

We argue by contradiction. If the statement were not true, there would exist a vector $c=(c_{1},c_{2},c_{3})\neq 0$ such that $Lc=0$, or equivalently
\begin{eqnarray}
\sum_{j=1}^{3}\bigg(\int_{\Sigma}\nu_{i}(x)\nu_{j}(x)d\sigma(x)\bigg)c_{j}=0.
\label{Lc=0}
\end{eqnarray}
Furthermore, writing $\nu_{i}$ as a linear combination of an orthonormal basis $\{e_{i}\}_{1\leq i\leq 3}$ of $span\{\nu_{i}\}_{1\leq i\leq 3}$, namely
\begin{eqnarray}\notag
\nu_{i}(x)=\sum_{k=1}^{3}\nu_{ik}e_{k}(x),
\end{eqnarray}
where
\begin{eqnarray}\notag
\nu_{ik}:=\int_{\Sigma}\nu_{i}(x)e_{k}(x)d\sigma(x),
\end{eqnarray}
we can see that, setting $a_{k}:=\sum_{j=1}^{3}\nu_{jk}c_{j}$, (\ref{Lc=0}) is equivalent to
\begin{eqnarray}\notag
0=\sum_{k=1}^{3}\bigg(\int_{\Sigma}\nu_{i}(x)e_{k}(x)d\sigma(x)\bigg)a_{k}=\int_{\Sigma}\nu_{i}(x)a(x)d\sigma(x)
\end{eqnarray}
with $a(x):=\sum_{k=1}^{3}a_{k}e_{k}(x)\in span\{e_{i}\}_{1\leq i\leq 3}=span\{\nu_{i}\}_{1\leq i\leq 3}$, so in particular $a\equiv 0$. On the other hand, $a_{k}=0$ for any $k$ is equivalent to
\begin{eqnarray}\notag
\int_{\Sigma}c(x)e_{k}(x)d\sigma(x)=0 &\text{ }\forall 1\leq k\leq 3,
\end{eqnarray}
with $c(x)=\sum_{j=1}^{3}c_{j}\nu_{j}(x)$. Thus $c\equiv 0$, that is $c_{j}=0$ for any $1\leq j\leq 3$, a contradiction.

\end{document}